\documentclass{amsart}
\usepackage{amsfonts}

\newtheorem{theorem}{Theorem}[section]
\newtheorem{lemma}[theorem]{Lemma}
\theoremstyle{definition}

\theoremstyle{remark}

\numberwithin{equation}{section}
\theoremstyle{plain}

\newtheorem{corollary}{Corollary}

\newtheorem{proposition}{Proposition}

\copyrightinfo{2001}{enter name of copyright holder}

\begin{document}
\title[Gradient Solitons]{On Gradient Ricci Solitons with Symmetry}
\author{Peter Petersen}
\address{520 Portola Plaza\\
Dept of Math UCLA\\
Los Angeles, CA 90095}
\email{petersen@math.ucla.edu}
\urladdr{http://www.math.ucla.edu/\symbol{126}petersen}
\thanks{}
\author{William Wylie}
\email{wylie@math.ucla.edu}
\urladdr{http://www.math.ucla.edu/\symbol{126}wylie}
\date{}
\subjclass{}
\keywords{}

\begin{abstract}
We study gradient Ricci solitons with maximal symmetry. First we show that
there are no non-trivial homogeneous gradient Ricci solitons. Thus the most
symmetry one can expect is an isometric cohomogeneity one group action. Many
examples of cohomogeneity one gradient solitons have been constructed.
However, we apply the main result in \cite{Petersen-Wylie} to show that
there are no noncompact cohomogeneity one shrinking gradient solitons with
nonnegative curvature.
\end{abstract}

\maketitle

\section{Introduction}

The goal of this paper is to study how symmetries can yield rigidity of a
gradient Ricci soliton together with weaker conditions than we used in \cite%
{Petersen-Wylie}. Recall that a Ricci soliton is a Riemannian metric
together with a vector field $\left( M,g,X\right) $ that satisfies%
\begin{equation*}
\mathrm{Ric}+\frac{1}{2}L_{X}g=\lambda g.
\end{equation*}%
It is called shrinking when $\lambda >0,$ steady when $\lambda =0$, and
expanding when $\lambda <0$. In case $X=\nabla f$ the equation can also be
written as%
\begin{equation*}
\mathrm{Ric}+\mathrm{Hess}f=\lambda g
\end{equation*}%
and is called a gradient (Ricci) soliton. A gradient soliton is rigid if it
is isometric to a quotient of $N\times \mathbb{R}^{k}$ where $N$ is an
Einstein manifold and $f=\frac{\lambda }{2}|x|^{2}$ on the Euclidean factor.
Throughout this paper we will also assume that our metrics have bounded
curvature. Shi's estimates for the Ricci flow then imply that all the
derivatives of curvature are also bounded (see Chapter 6 of \cite{Chow-Lu-Ni}%
).

First we show that all gradient solitons with maximal symmetry are rigid.

\begin{theorem}
\label{homogeneous} All homogeneous gradient Ricci solitons are rigid.
\end{theorem}

This is in sharp contrast to the more general Ricci solitons that exist on
many Lie groups and other homogeneous spaces see \cite{Baird-Danielo,  Lauret2001, Lott}. It also
shows that the maximal amount of symmetry we can expect on a nontrivial
gradient soliton is a cohomogeneity 1 action that leaves $f$ invariant.
Particular cases, such as the rotationally symmetric case on $\mathbb{R}^{n}$
and the $U(n)$ invariant case on certain K\"{a}hler manifolds, have been
studied extensively and many interesting examples have been found, see e.g. 
\cite{Cao1996, Cao1997, Dancer-Wang, Feldman-Ilmanen-Knopf2003,Ivey1994, Koiso1990,
Kotschwar, Wang-Zhu2004, Yang}. In particular, Kotschwar \cite{Kotschwar} has
shown that the only  rotationally symmetric shrinking gradient soliton metrics on $%
S^{n},$ $\mathbb{R}^{n},$ and $S^{n-1}\times \mathbb{R}$ are the rigid ones and  Feldman, Ilmanen, and Knopf \cite{Feldman-Ilmanen-Knopf2003} have
 proven that the only $U(n)$ invariant shrinking soliton on $\mathbb{C}%
^{n}$ is the flat metric. No curvature assumption is required for these results.  On the other hand, there  are non-rigid complete noncompact $U(n)$ invariant
gradient shrinking solitons \cite{Dancer-Wang, Feldman-Ilmanen-Knopf2003, Yang}. These examples  show that some other assumption is necessary in general to prove rigidity.  Here we show that nonnegative curvature suffices. 

\begin{theorem}
\label{cohom1} All complete noncompact shrinking gradient solitons of
cohomogeneity 1 with nonnegative Ricci curvature and $\mathrm{sec}\left(
E,\nabla f\right) \geq 0$ are rigid.
\end{theorem}

Recently  Naber \cite{Naber},  building on work of Ni and Wallach \cite{Ni-WallachII}, has shown that every 4-dimensional complete shrinking soliton with nonnegative curvature operator is rigid. (This was proven in dimensions 2 and 3 by  Hamilton \cite{Hamilton1988} and Perelman \cite{PerelmanII} respectively.)  Theorem \ref{cohom1} offers further evidence this result extends to higher dimensions.  In fact, in the proof all we use about  cohomogeneity one is a much weaker condition on $f$ we call rectifiability which we will discuss in section 3.  (Also recall that the work of  B\"{o}hm and Wilking \cite{Bohm-Wilking}  implies  that every compact shrinking gradient Ricci soliton with nonnegative curvature operator is a quotient of the round sphere.)  For other recent results concerning the classification of gradient shrinking solitons see \cite{Eminenti-LaNave-Mantegazza, Ni, Ni-WallachI, Petersen-WylieIII, Weber, Wylie}.  

The famous Bryant soliton (see \cite{Ivey1994}) and the examples in \cite{Cao1997} show that there are non-rigid rotationally symmetric steady and expanding gradient solitons with positive curvature operator.

\section{Killing Fields on Gradient Solitons}

In this section we establish a splitting theorem involving Killing fields on
a gradient soliton which leads to Theorem \ref{homogeneous}. The main
observation is the following.

\begin{proposition}
If $X$ is a Killing field on a gradient soliton, then $\nabla D_{X}f$ is
parallel. Moreover, if $\lambda \neq 0$ and $\nabla D_{X}f=0$ then also $%
D_{X}f=0.$
\end{proposition}

\begin{proof}
We have that $L_{X}g=0,$ thus $L_{X}\mathrm{Ric}=0$ and hence%
\begin{eqnarray*}
0 &=&L_{X}\mathrm{Hess}f \\
&=&\mathrm{Hess}L_{X}f \\
&=&\mathrm{Hess}D_{X}f
\end{eqnarray*}%
this proves the first claim.

Next note that if $\nabla D_{X}f=0,$ then $D_{X}f$ is constant. Thus $f\circ
\gamma _{X}\left( t\right) :\mathbb{R\rightarrow R}$ is onto if $\gamma _{X}$
is an integral curve for $X$ and $D_{X}f$ doesn't vanish.

On the other hand recall that the soliton equation implies that 
\begin{equation*}
\mathrm{scal}+\left\vert \nabla f\right\vert ^{2}-2\lambda f=\mathrm{const}
\end{equation*}%
So if the scalar curvature is bounded we see that $f$ must either be bounded from
below or above and hence $D_{X}f=0.$
\end{proof}

This shows that either $D_Xf=0$ or the metric splits off a Euclidean factor.
One might worry that the soliton structure may not also split, however the
next lemma shows this is not an issue.

\begin{lemma}
If a gradient soliton splits $\left( M,g\right) =\left( M_{1}\times
M_{2},g_{1}+g_{2}\right) $ as a Riemannian product, then $f\left(
x_{1},x_{2}\right) =f_{1}\left( x_{1}\right) +f_{2}\left( x_{2}\right) $
also splits in such a way that each $\left( M_{i},g_{i},f_{i}\right) $ is a
soliton 
\begin{equation*}
\mathrm{Ric}_{g_{i}}+\mathrm{Hess}f_{i}=\lambda g_{i}
\end{equation*}
\end{lemma}

\begin{proof}
Use the $\left( 1,1\right) $ version of the soliton equation%
\begin{equation*}
\mathrm{Ric}+\nabla \nabla f=\lambda I
\end{equation*}%
to see that the operator $E\rightarrow \nabla _{E}\nabla f$ preserves the
manifold splitting as the Ricci curvature preserves the splitting. This can
be used to first split the gradient $\nabla f.$ To see how, use local
coordinates $x^{j}$ such $x^{1},...,x^{m}$ are coordinates on $M_{1}$ and $%
x^{m+1},...,x^{n}$ coordinates on $M_{2}$. The splitting of the metric then
implies that%
\begin{equation*}
\nabla _{\partial _{i}}\partial _{j}=0
\end{equation*}%
if $i\leq m$ and $j\geq m+1$ or $i\geq m=1$ and $j\leq m.$ If we write $%
\nabla f=\alpha ^{j}\partial _{j},$ then%
\begin{eqnarray*}
\nabla _{\partial _{i}}\nabla f &=&\nabla _{\partial _{i}}\alpha
^{j}\partial _{j} \\
&=&\left( \partial _{i}\alpha ^{j}\right) \partial _{j}+\alpha ^{j}\nabla
_{\partial _{i}}\partial _{j}.
\end{eqnarray*}%
If we assume that $i\leq m$ then%
\begin{eqnarray*}
\nabla _{\partial _{i}}\nabla f &\in &TM_{1}, \\
\alpha ^{j}\nabla _{\partial _{i}}\partial _{j} &\in &TM_{1}
\end{eqnarray*}%
showing that $\partial _{i}\alpha ^{j}=0$ for $j\geq m+1.$ Similarly $%
\partial _{i}\alpha ^{j}=0$ when $i\geq m+1$ and $j\leq m.$ This shows that 
\begin{equation*}
\nabla f=X_{1}+X_{2}
\end{equation*}%
where $X_{i}$ are vector fields on $M_{i}.$ We then see that 
\begin{equation*}
X_{i}=\nabla f_{i}
\end{equation*}%
where 
\begin{eqnarray*}
f_{1}\left( x_{1}\right) &=&f\left( x_{1},q\right) , \\
f_{2}\left( x_{2}\right) &=&f\left( p,x_{2}\right) -f\left( p,q\right)
\end{eqnarray*}%
for some fixed point $\left( p,q\right) \in M_{1}\times M_{2}.$
\end{proof}

Note that the splitting of the metric implies%
\begin{equation*}
R\left(\cdot ,\nabla f\right) \nabla f=R_{1}\left(\cdot ,\nabla f_{1}\right)
\nabla f_{1}+R_{2}\left(\cdot ,\nabla f_{2}\right) \nabla f_{2}
\end{equation*}%
So if, say, $M_{2}$ is flat then the radial curvatures of $M$ and $M_{1}$
are the same.

This implies the reduction result alluded to above.

\begin{corollary}
\label{reduction} If $X$ is a Killing field on a gradient soliton, then
either $D_{X}f=0$ or we have an isometric splitting $M=N\times \mathbb{R}$
where $N$ is a gradient soliton with the same radial curvatures as $M.$
\end{corollary}

Intuitively, Corollary \ref{reduction} says that if the metric of a gradient
soliton has some symmetry, then the only way $f$ can break the symmetry is
by splitting off a Gaussian factor. With this fact we can prove the result
for homogeneous solitons.

\begin{theorem}
All homogeneous gradient solitons are rigid.
\end{theorem}

\begin{proof}
In case the soliton is steady this is a consequence of the scalar curvature
being constant and hence $M$ is Ricci flat.

When the soliton is expanding or shrinking split $M=N\times \mathbb{R}^{k}$
such that $N$ doesn't have any flat de Rham factors. If $G$ acts
transitively on $M$ it also acts transitively on each of the two factors as
isometries preserve the flat de Rham factor.

The previous lemma and corollary now tell us that all Killing fields on $N$
must leave $f_{1}$ invariant. Thus $N$ can't be homogeneous unless $f_{1}$
is trivial.
\end{proof}

\section{Rectifiability}

In this section we prove the result for cohomogeneity one and more general
rectifiable gradient solitons.

We say that a function $u$ is \emph{rectifiable} if it can be written as $%
u=h\left( r\right) $ where $r$ is a distance function. It is easy to check
that a function is rectifiable if and only if its gradient $\nabla u$ has
constant length on the level sets of $u.$ We will say that a gradient soliton $%
(M,g,f)$ is rectifiable if the function $f$ is rectifiable on $(M,g)$.

It is easy to see that a gradient soliton with a cohomogeneity 1 group
action that leaves $f$ invariant is rectifiable. Assume that $G$ is such a
isometric group action. This gives us a distance function 
\begin{equation*}
r:M\rightarrow M/G\subset \mathbb{R}
\end{equation*}%
(locally if $G$ is noncompact) and $f=h\left( r\right) $ as $f$ is constant
on the orbits of the action. Similarly the scalar curvature is also
rectifiable with respect to $r.$

We note the following interesting properties of rectifiable solitons.

\begin{proposition}
If $(M,g,f)$ is a rectifiable gradient soliton with $f=h(r)$ then $\mathrm{%
scal}$, $\Delta f$, and $\Delta r$ are also rectifiable. In particular, $%
\mathrm{Ric}(\nabla f, \nabla f)=0$ if and only if $(M,g)$ has constant
scalar curvature.
\end{proposition}

\begin{proof}
If $f$ is rectifiable, then $|\nabla f|$ is also rectifiable so the equation 
\begin{equation*}
\mathrm{scal}+|\nabla f|^{2}-2\lambda f=\mathrm{const}
\end{equation*}
implies that the scalar curvature is rectifiable.

Tracing the soliton equation then gives 
\begin{equation*}
\mathrm{scal}=\lambda n-\Delta f,
\end{equation*}
so $\Delta f$ is rectifiable. Since $f$ is rectifiable we can write 
\begin{equation*}
\Delta f=h^{\prime \prime }(r)+h^{\prime }(r)\Delta r
\end{equation*}%
so $\Delta f$ rectifiable implies that $\Delta r$ is also rectifiable.

Now since $\mathrm{scal}$ and $f$ are rectifiable $\nabla \mathrm{scal} = 
\mathrm{Ric}(\nabla f)$ is proportional to $\nabla f$, proving the last
statement.
\end{proof}

The main result from \cite{Petersen-Wylie} now shows that a rectifiable
gradient soliton is rigid if and only if it is radially flat. We note that,
in the case of cohomogeneity one, radial flatness, even without the soliton
equation, is already quite restrictive.

\begin{theorem}
A radially flat cohomogeneity 1 space coming from a compact action is a flat
bundle.
\end{theorem}

\begin{proof}
Let $r:M\rightarrow \mathbb{R}$ be the distance function coming from the
quotient $M\rightarrow M/G.$ It is smooth except at the singular orbits. The
singular orbits correspond to the minimum and/or maximum of $r$ if they
exist.

Let $S_{r}=\nabla \nabla r,$ then%
\begin{equation*}
\nabla _{\nabla r}S_{r}+S_{r}^{2}=0
\end{equation*}%
This means that $S_{r}$ is completely determined by the singular orbits
where $S_{r}\rightarrow 0$ on vectors tangent to the singular orbit and $%
S_{r}\rightarrow \infty $ on vectors normal to the singular orbit and
perpendicular to $\nabla r.$

If there are no singular orbits, then $S_{r}=0$ is the only possibility as
all other solutions blow up in finite time going forwards or backwards. Thus
the space splits.

If $r$ has a minimum set, then solutions that start out being zero stay
zero, while the other solutions that start out being $\infty $ decay to
zero. As they never become zero the space is noncompact. We see that the
space must then be a flat bundle $N\times _{\Gamma }\mathbb{R}^{k}$ where $%
N/\Gamma $ is the singular orbit.
\end{proof}

We now turn our attention to proving rigidity for rectifiable shrinking
solitons with nonnegative radial curvature.

\begin{proposition}
Let $\left( M,g\right) $ be a Riemannian manifold and $r:M\rightarrow
\lbrack 0,\infty )$ a proper distance function that is smooth outside a
compact set. If $\mathrm{sec}\left( E,\nabla r\right) \geq 0,$ then $r$ is
convex outside a compact set.
\end{proposition}

\begin{proof}
Define $S_{r}=\nabla \nabla r$ and use that it solves the equation%
\begin{equation*}
\nabla _{\nabla r}S_{r}=-S_{r}^{2}-R\left( \cdot ,\nabla r\right) \nabla r.
\end{equation*}%
As $E\rightarrow R\left( E,\nabla r\right) \nabla r$ is assumed to be
nonnegative we see that if $S_{r}$ has a negative eigenvalue somewhere, then
it will go to $-\infty $ before $r$ reaches infinity. This contradicts that $%
r$ is smooth.
\end{proof}

\begin{lemma}
Let $\left( M,g,f\right) $ be a noncompact nontrivial shrinking gradient
soliton with rectifiable and proper $f$. If the radial curvatures, $\mathrm{%
sec}\left( E,\nabla f\right) $ are nonnegative, then $f$ is convex at
infinity.
\end{lemma}

\begin{proof}
Since $f$ is rectifiable: $f=h\left( r\right) ,$ where $r:M\rightarrow
\lbrack 0,\infty )$ is a distance function that is smooth outside a compact
set. Since $f$ and $r$ have proportional gradients, $\nabla f=h^{\prime
}\nabla r,$ our curvature assumption guarantees that $r$ is convex at
infinity.

First note that the equation%
\begin{equation*}
\mathrm{scal}+\left\vert \nabla f\right\vert ^{2}-2\lambda f=\mathrm{const}
\end{equation*}%
shows that $\left\vert \nabla f\right\vert \rightarrow \infty $ as $%
\mathrm{scal}$ is bounded and $f$ is proper, i.e., $\left\vert f\right\vert
\rightarrow \infty .$ In particular $h^{\prime }>0$ outside a compact set.

Define $S_{f}=\nabla \nabla f$ and $S_{r}=\nabla \nabla r,$ they are related
by%
\begin{eqnarray*}
S_{f} &=&\nabla \nabla f \\
&=&h^{\prime \prime }dr\otimes \nabla r+h^{\prime }\nabla \nabla r \\
&=&h^{\prime \prime }dr\otimes \nabla r+h^{\prime }S_{r}
\end{eqnarray*}%
The soliton equation shows that%
\begin{equation*}
\mathrm{Ric}\left( \nabla r,\nabla r\right) +h^{\prime \prime }=\lambda
\end{equation*}%
Since $\mathrm{Ric}\left( \nabla r,\nabla r\right) $ is nonnegative this
shows that $h^{\prime \prime }\leq \lambda $.

Next we claim that $\mathrm{Ric}\left( \nabla r\right) \rightarrow 0$ as $%
r\rightarrow \infty .$ This follows from the formula%
\begin{eqnarray*}
\mathrm{Ric}\left( \nabla r\right)  &=&\pm \frac{\mathrm{Ric}\left( \nabla
f\right) }{\left\vert \nabla f\right\vert } \\
&=&\pm \frac{1}{2}\frac{\nabla \mathrm{scal}}{\left\vert \nabla f\right\vert 
}
\end{eqnarray*}%
where we note that $\nabla \mathrm{scal}$ is bounded and $\left\vert \nabla
f\right\vert \rightarrow \infty $ at infinity.

Thus%
\begin{equation*}
S_{f}\left( \nabla r\right) =h^{\prime \prime }\nabla r\sim \lambda \nabla r
\end{equation*}%
at infinity. This proves that outside some large compact set $h^{\prime
\prime }\geq \lambda /2$ and $h^{\prime }>0.$ Thus $f$ is convex outside a
compact set.
\end{proof}

\begin{theorem}
A complete, noncompact, rectifiable, shrinking gradient soliton with
nonnegative radial sectional curvature, and nonnegative Ricci curvature is
rigid.
\end{theorem}

\begin{proof}
Let $f=h(r)$. Since we have a shrinking gradient Ricci soliton with bounded nonnegative curvature $f$ is proper \cite{PerelmanII}.  Therefore, the previous lemmas show that $f$ and $r$ are proper and convex outside a compact set. This
implies that $\mathrm{Ric}\leq \lambda g$ outside a compact set. Define $%
\Delta _{f}=\Delta -D_{\nabla f}$ to be the $f$-Laplacian, then (see \cite%
{Petersen-Wylie}) 
\begin{equation*}
\Delta _{f}\mathrm{scal}=\mathrm{tr}\left( \mathrm{Ric}\circ (\lambda I-%
\mathrm{Ric})\right)
\end{equation*}%
So $\mathrm{Ric}\leq \lambda g$ outside a compact set implies%
\begin{equation*}
\Delta _{f}\mathrm{scal}\geq 0
\end{equation*}%
outside a set $\Omega _{R}=\left\{ x\in M:r\leq R\right\} .$ We also know
that $\mathrm{scal}$ is increasing along gradient curves for $\nabla f$ as%
\begin{equation*}
D_{\nabla f}\mathrm{scal}=2\mathrm{Ric}\left( \nabla f,\nabla f\right) \geq
0.
\end{equation*}%
If 
\begin{equation*}
s_{R}=\min_{p\in \partial \Omega _{R}}\mathrm{scal}_{p}
\end{equation*}%
then the function%
\begin{equation*}
u=\max \left\{ \mathrm{scal},s_{R}\right\}
\end{equation*}%
satisfies%
\begin{equation*}
\Delta _{f}u\geq 0
\end{equation*}%
From Theorem  4.2 in \cite{Petersen-WylieIII} it follows that $u$ is constant (also
see \cite{Wei-Wylie}). This shows that $\mathrm{scal}=s_{R}$ on $%
M-\Omega _{R}.$ Since $\left( M,g\right) $ is analytic (see \cite{Bando1987}%
) the scalar curvature is constant on all of $M.$ This in turn shows that $%
\mathrm{Ric}\left( \nabla f,\nabla f\right) =0$ everywhere and hence $%
\mathrm{sec}(E,\nabla f)\geq 0$ implies that $\left( M,g\right) $ is
radially flat. The main theorem from \cite{Petersen-Wylie} then shows that $%
M $ is rigid.
\end{proof}


\begin{thebibliography}{99}

\bibitem{Baird-Danielo} Paul Baird and Laurent Danielo. \newblock %
Three-dimensional Ricci solitons which project to surfaces. 
\newblock {\em J. Reine Angew. Math.} 6008(2007), 65-91.

\bibitem{Bando1987} Shigetoshi Bando. \newblock Real Analyticity of
Solutions of Hamilton's equation. \newblock Math Z. (195): 93-97, 1987.

\bibitem{Bohm-Wilking} Christoph B\"{o}hm and Burkhard Wilking.  \newblock Manifolds with positive curvature operators are space forms. \newblock To appear in Ann. of Math.  \newblock arXiv:math/0606187.



\bibitem{Cao1996} Huai-Dong Cao. \newblock Existence of gradient {K}\"{a}%
hler-{R}icci solitons. \newblock In \emph{Elliptic and parabolic methods in
geometry (Minneapolis, MN, 1994)}, pages 1--16. A K Peters, Wellesley, MA,
1996.

\bibitem{Cao1997} Huai-Dong Cao. \newblock Limits of solutions to the K\"{a}%
hler-Ricci flow. \newblock {\em J. Differential Geom.}, 45(2):257-272 1997.



\bibitem{Chow-Lu-Ni} Bennet Chow, Peng Lu, and Lei Ni. \newblock Hamilton's
Ricci flow. \newblock {\em Graduate studies in Mathematics}, AMS,
Providence, RI, 2006.

\bibitem{Dancer-Wang} Andrew Dancer and Mckenzie Wang. \newblock On Ricci Solitons of cohomogeneity one. \newblock arXiv:math/0802.0759.

\bibitem{Eminenti-LaNave-Mantegazza} Manolo Eminenti, Gabriele La Nave, and
Carlo Mantegazza. \newblock Ricci Solitons - the Equation Point of View. %
\newblock arXiv:math.DG/0607546v2. 


\bibitem{Feldman-Ilmanen-Knopf2003} Mikhail Feldman, Tom Ilmanen, and Dan
Knopf. \newblock Rotationally symmetric shrinking and expanding gradient {K}%
\"{a}hler-{R}icci solitons. \newblock {\em J. Differential Geom.},
65(2):169--209, 2003.

\bibitem{Hamilton1988} Richard Hamilton. \newblock The Ricci flow on
Surfaces. \newblock In \emph{Mathematics and General Relativity.} volume 71
of Contemporary Mathematics pages 237-262. AMS, Providence, RI, 1988.


\bibitem{Ivey1994} Thomas Ivey. \newblock New Examples of complete Ricci
solitons. \newblock{ \em Proc. Amer. Math Soc.} 122(1): 241-245, 1994.

\bibitem{Koiso1990} Norihito Koiso. \newblock On rotationally symmetric {H}%
amilton's equation for {K}\"ahler-{E}instein metrics. \newblock In \emph{%
Recent topics in differential and analytic geometry}, volume~18 of \emph{%
Adv. Stud. Pure Math.}, pages 327--337. Academic Press, Boston, MA, 1990.

\bibitem{Kotschwar} Brett Kotschwar. \newblock On rotationally invariant
shrinking gradient Ricci solitons. \newblock arXiv:math/0702597.

\bibitem{Lauret2001} Jorge Lauret. \newblock Ricci Soliton homogeneous
nilmanifolds. \newblock{\em Math. Ann.}, 319:715-733, 2001.

\bibitem{Lott} John Lott. \newblock On the long-time behavior of type-III
Ricci flow solutions. \newblock {\em Math. Ann.} 339(2007), no. 3, 627-666.

\bibitem{Naber} Aaron Naber. \newblock Noncompact shrinking $4$-solitons
with nonnegative curvature. \newblock arXiv:math.DG/0710.5579.

\bibitem{Ni} Lei Ni. \newblock Ancient Solutions to K\"ahler Ricci flow. %
\newblock {\em Math. Research Letters.}, 12: 633-654, 2005.

\bibitem{Ni-WallachI} Lei Ni and Nolan Wallach. \newblock On a classification
of the gradient shrinking solitons. \newblock arXiv:math.DG/0710.3194.

\bibitem{Ni-WallachII} Lei Ni and Nolan Wallach. \newblock On 4-dimensional gradient shrinking solitons. \newblock arXiv:math.DG/0710.3195.




\bibitem{PerelmanII} G.~Ya. Perelman. \newblock Ricci flow with surgery on
three manifolds. \newblock arXiv: math.DG/0303109.


\bibitem{Petersen-Wylie} Peter Petersen and William Wylie. \newblock Rigidity of
gradient Ricci solitons. \newblock  arXiv:math.DG/0710.3174.

\bibitem{Petersen-WylieIII} Peter Petersen and William Wylie. \newblock On the classification of gradient Ricci solitons. \newblock arXiv:math.DG/0712.1298.


\bibitem{Weber} Brian Weber. \newblock Convergence of compact Ricci solitons. \newblock arXiv:math.{DG}/0804.1158.

\bibitem{Wei-Wylie} Guofang Wei and William Wylie. \newblock Comparison
Geometry for the Bakry-Emery Ricci tensor. \newblock  arXiv:math.{DG}%
/0706.1120.

\bibitem{Wylie} William Wylie. \newblock Complete shrinking {R}icci solitons
have finite fundamental group. \newblock {\em Proc. American Math. Sci.}
136(2008), no. 5, 1803-1806.

\bibitem{Wang-Zhu2004} X.J. Wang and X.H. Zhu. \newblock K\"ahler-Ricci
solitons on toric manifolds with positive first Chern class. 
\newblock {\em
Adv. Math.}, 188(1):87-103, 2004.

\bibitem{Yang} Bo Yang. \newblock A characterization of Koiso's typed solitons. \newblock arXiv:math.{DG}/0802.0300. 

\end{thebibliography}
\end{document}